\documentclass[a4paper,11pt,reqno]{amsart}

\usepackage{latexsym,dsfont,amssymb,amsmath,amsthm}

\chardef\bslash=`\\ 

\numberwithin{equation}{section}

\newtheorem{theorem}{Theorem}[section]
\newtheorem{corollary}[theorem]{Corollary}
\newtheorem{lemma}[theorem]{Lemma}
\newtheorem{proposition}[theorem]{Proposition}

\theoremstyle{remark}
\newtheorem{remark}[theorem]{Remark}

\theoremstyle{definition}

\newcommand\bp{\begin{proof}}
\newcommand\ep{\end{proof}}

\newcommand\3[1]{{\mathds #1}}
\newcommand\inv{^{-1}}

\newcommand\pp{\mathfrak p}
\newcommand\qq{\mathfrak q}

\newcommand{\N}{\mathbb N}
\newcommand{\Z}{\mathbb Z}

\newcommand{\HH}{\mathbb H}
\newcommand{\Q}{\mathbb Q}
\newcommand{\R}{\mathbb R}

\newcommand\T{\mathbb T}
\newcommand\A{\mathbb{A}}
\newcommand\af{\mathbb{A}_f}
\newcommand\ak{{\mathbb A}_K}
\newcommand\akf{{\mathbb A}_{K,f}}
\newcommand\jkf{\akf^*}
\newcommand\OO{{\mathcal O}}
\newcommand\ohs{{\hat{\OO}^*}}
\newcommand\gal{\mathcal G}
\newcommand\kab{K^{ab}}

\newcommand{\hh}{\mathcal H}
\newcommand{\rr}{\mathcal R}
\newcommand{\primes}{\mathcal  P}

\newcommand\eps{\varepsilon}

\newcommand\bpmatrix{\begin{pmatrix}}
\newcommand\epmatrix{\end{pmatrix}}

\newcommand{\diag}[2]{\left(\begin{matrix}#1&0\\
0&#2\end{matrix}\right)}

\newcommand\mz{{\operatorname{Mat}_2}(\Z)}

\newcommand\ma{{\operatorname{Mat}_2}(\af)}
\newcommand\mtwo{\operatorname{Mat}_2}

\newcommand\glq{{\operatorname{GL}^+_2}(\Q)}

\newcommand\slz{{\operatorname{SL}_2}(\Z)}

\newcommand\gl{{\operatorname{GL}_2}}
\newcommand\glp{{\operatorname{GL}_2^+}}
\newcommand\sltwo{{\operatorname{SL}_2}}
\newcommand\gln{{\operatorname{GL}_n}}

\newcommand\pgl{{\operatorname{PGL}^+_2}(\R)}
\newcommand\pgll{{\operatorname{PGL}_2}(\R)}

\newcommand\pso{{\operatorname{PSO}_2}(\R)}

\newcommand{\hecke}[2]{\mathcal \hh({#1}, {#2})}

\newcommand{\red}[1]{C^*_r(#1)}

\newcommand{\fib}[3]{#2\backslash #1\times_{#2}#3}
\newcommand{\fibb}[3]{#2\backslash #1\boxtimes_{#2}#3}

\newcommand\enu[1]{\smallskip\newline\makebox[5mm][l]{\rm(#1)}}

\begin{document}

\title[Bost-Connes type systems]
{Von Neumann algebras arising from Bost-Connes type systems}

\author[S. Neshveyev]{Sergey Neshveyev}
\address{Department of Mathematics, University of Oslo,
P.O. Box 1053 Blindern, N-0316 Oslo, Norway\\ Institut de
Math\'ematiques de Jussieu, Universit\'e Paris 7 Denis Dedirot, 175 rue
du Chevaleret, 75013 Paris, France}

\email{sergeyn@math.uio.no}

\thanks{Supported by the Research Council of Norway.}

\begin{abstract}
We show that the KMS$_\beta$-states of Bost-Connes type systems for
number fields in the region $0<\beta\le1$, as well as of the
Connes-Marcolli $\gl$-system for $1<\beta\le2$, have type III$_1$. This
is equivalent to ergodicity of various actions on adelic spaces. For
example, the case $\beta=2$ of the $\gl$-system corresponds to
ergodicity of the action of $\gl(\Q)$ on $\mtwo(\A)$ with its Haar
measure.
\end{abstract}

\date{July 9, 2009}

\maketitle

\bigskip

\section*{Introduction}

The Bost-Connes system~\cite{bos-con} is a C$^*$-dynamical system such
that for inverse temperatures $\beta>1$ the extremal KMS$_\beta$-states
carry a free transitive action of the Galois group of the maximal
abelian extension of $\Q$, have type I and partition function
$\zeta(\beta)$, while for every $\beta\in(0,1]$ there exists a unique
KMS$_\beta$-state of type III$_1$. The uniqueness and the type of the
KMS$_\beta$-states in the critical interval~$(0,1]$ are the most
difficult parts of the analysis of the system. The result is equivalent
to ergodicity of certain measures on the space $\A$ of adeles with
respect to the action of $\Q^*$. In particular, the case $\beta=1$
corresponds to a Haar measure $\mu_1$.
To see why ergodicity of $\mu_1$ is nontrivial, observe that as $\Q$ is
discrete in $\A$, the orbit of any point in $\A^*$ is discrete in $\A$,
while the ergodicity implies that almost every orbit in $\A$ is dense.
There is of course no contradiction since $\A^*$ is a subset of $\A$ of
measure zero, but one does realize that it is difficult to immediately
see a single dense orbit in $\A$.

Recently the construction of Bost and Connes has been generalized first
to imaginary quadratic fields~\cite{cmn1} and then to arbitrary number
fields~\cite{ha-pa}. The crucial step of imaginary quadratic fields was
achieved by introducing a universal system of quadratic fields, the so
called $\gl$-system of Connes and Marcolli~\cite{con-mar}.
In~\cite{LLNlat} and~\cite{LLNcm} we analyzed these systems in the
critical intervals $0<\beta\le1$ for number fields and $1<\beta\le2$
for the $\gl$-system, and showed that there exist unique
KMS$_\beta$-states. The aim of the present paper is to prove that these
KMS$_\beta$-states have type III$_1$. For number fields the proof is
similar to the one for the original Bost-Connes
system~\cite{bos-con,nes}. The interesting case is that of the
$\gl$-system. It amounts to proving that the action of $\gl(\Q)$ on
$\pgll\times\ma$ has type~III$_1$ with respect to certain
product-measures. After passing to the quotient space
$\gl(\Z)\backslash(\pgll\times\ma)/\gl(\hat\Z)$, we essentially use an
argument showing that a nonzero element of the asymptotic ratio set is
contained in the ratio set, written in terms of $L^2$-spaces rather
than measure spaces and using a representation of the Hecke algebra
$\hecke{\gl(\Q)}{\gl(\Z)}$ instead of group actions. An additional
difficulty is that we do not have a product decomposition of the
representation of the Hecke algebra. In other words, the Hecke
operators defined by elements of $\gl(\Z[p^{-1}])$ act nontrivially on
$\gl(\Z)\backslash(\pgll\times\prod_{q\ne p}\mtwo(\Z_q))/\gl(\hat \Z)$.
As a side remark, a similar problem would not arise for the finite part
of the $\gl$-system~\cite{LLNfcm}. What saves the day for the full
system is that this action is mixing on large subsets, which is a
consequence of a variant of equidistribution of Hecke
points~\cite{cou}.

Apart from some trivial cases when there is a subgroup of measure
preserving transformations acting ergodically (for example, for
$\gl(\Q)$ acting on $\R^2$), computations of ratio sets are usually
quite hard, see e.g. \cite{bow,hs,ino}. A large class of type~III$_1$
actions can be obtained as follows~\cite{zim}. Let $G$ be a connected
non-compact simple Lie group with finite center, $\Gamma\subset G$ a
lattice and $P\subset G$ a parabolic subgroup. Then the action of
$\Gamma$ on $G/P$ has type III$_1$. This is proved by identifying the
underlying space of the associated flow with the measure-theoretic
quotient $\Gamma\backslash G/P_0$, where $P_0$ is the kernel of the
modular function of $P$, and using that the action of $P_0$ on
$\Gamma\backslash G$ is mixing by Howe-Moore's theorem~\cite{zim}. With
these examples in mind, it seems only natural that to compute the type
of the states of the $\gl$-system a form of adelic mixing is needed.

\medskip

\noindent{\bf Acknowledgement.} It is my pleasure to thank Hee Oh for
her help with equidistribution of Hecke points.

\bigskip

\section{Actions of type III$_1$}

Assume a countable group $G$ acts ergodically on a measure space
$(X,\mu)$. The ratio set of the action~\cite{kr} consists of all
numbers $\lambda\ge0$ such that for any $\eps>0$ and any subset
$A\subset X$ of positive measure there exists $g\in G$ such that
$$
\mu\left(\left\{x\in gA\cap A:\left|\frac{dg\mu}{d\mu}(x)
-\lambda\right|<\eps\right\}\right)>0,
$$
where $g\mu$ is the measure defined by $g\mu(Z)=\mu(g^{-1}Z)$. The
ratio set depends only on the orbit equivalence relation
$\rr=\{(x,gx)\mid x\in X,\ g\in G\}\subset X\times X$. We will denote
it by $r(\rr)$. The set $r(\rr)\setminus\{0\}$ is a closed subgroup of
$\R^*_+$. The action is said to be of type III$_1$ if this subgroup
coincides with the whole group $\R^*_+$.

Denote by $\lambda_\infty$ the Lebesgue measure on $\R$. We have two
commuting actions of $\R$ and $G$ on $(\R_+\times
X,\lambda_\infty\times\mu)$,
$$
g(t,x)=\left(\frac{dg\mu}{d\mu}(gx)t,gx\right)\ \ \text{for}\ \ g\in G,\ \
s(t,x)=(e^{-s}t,x)\ \ \text{for}\ \ s\in\R.
$$
The flow of weights~\cite{CT} of the von Neumann algebra $W^*(\rr)$ is
the flow induced by the above action of $\R$ on the measure-theoretic
quotient of $(\R_+\times X,\lambda_\infty\times\mu)$ by the action of
$G$. The original action of $G$ on $(X,\mu)$ has type III$_1$ if and
only if the flow of weights is trivial, that is, the action of $G$ on
$(\R_+\times X,\lambda_\infty\times\mu)$ is ergodic.

\medskip

Let $\{(X_n,\mu_n)\}_{n=1}^\infty$ be a sequence of at most countable
probability spaces. Put $(X,\mu)=\prod_n(X_n,\mu_n)$, and define an
equivalence relation $\rr$ on $X$ by
$$
x\sim y\ \ \text{if}\ \ x_n=y_n\ \ \text{for all} \ \ n
\ \ \text{large enough}.
$$
For a finite subset $I\subset\N$ and $a\in\prod_{n\in I}X_n$ put
$$
Z(a)=\{x\in X\mid x_n=a_n\ \text{for}\ n\in I\}.
$$
The asymptotic ratio set $r_\infty(\rr)$ consists by
definition~\cite{AW} of all numbers $\lambda\ge0$ such that for any
$\eps>0$ there exist a sequence $\{I_n\}^\infty_{n=1}$ of mutually
disjoint finite subsets of $\N$, disjoint subsets
$K_n,L_n\subset\prod_{k\in I_n}X_k$ and bijections $\varphi_n\colon
K_n\to L_n$ such that
$$
\left|\frac{\mu(Z(\varphi_n(a)))}{\mu(Z(a))}-\lambda\right|<\eps\ \text{for all}\
a\in K_n\ \text{and}\ n\ge1,\ \ \text{and}\ \ \sum_{n=1}^\infty\sum_{a\in K_n}\mu(Z(a))=\infty.
$$
It is known that $r_\infty(\rr)\setminus\{0\}=r(\rr)\setminus\{0\}$. We
will only need the rather obvious inclusion $\subset$.

\bigskip

\section{Bost-Connes type systems for number fields}

Suppose $K$ is an algebraic number field with subring of integers
$\OO$. Denote by $V_K$ the set of places of $K$, and by $V_{K,f}\subset
V_K$ the subset of finite places. For $v\in V_K$ denote by~$K_v$ the
corresponding completion of $K$. If $v$ is finite, let $\OO_v$ be the
closure of $\OO$ in $K_v$. Denote also by
$K_\infty=\prod_{v|\infty}K_v$ the completion of $K$ at all infinite
places. The adele ring $\ak$ is the restricted product of the rings
$K_v$ with respect to $\OO_v\subset K_v$, $v\in V_K$. When the product
is restricted to $v\in V_{K,f}$, we get the ring $\akf$ of finite
adeles. The ring of finite integral adeles is $\hat\OO=\prod_{v\in
V_{K,f}}\OO_v\subset \akf$. We identify $\akf^*$ with the subgroup of
$\ak^*$ consisting of elements with coordinates $1$ for all infinite
places.

Consider the topological space $\gal(\kab / K)\times \akf$, where
$\gal(\kab/K)$ is the Galois group of the maximal abelian extension of
$K$. On this space there is an action of the group $\jkf$ of finite
ideles, via the Artin map $s\colon \ak^* \to \gal(\kab/K)$ on the first
component and via multiplication on the second component:
\[
j(\gamma, m) = (\gamma s(j)\inv , jm)  \ \ \hbox{for}\ \ j\in \jkf,
\ \ \gamma\in \gal(\kab/K), \ \ m\in \akf.
\]
Consider the quotient space $\gal(\kab / K)\times_{\ohs} \akf$ by the
action of $\ohs \subset \ak^*$. On this space we have a quotient action
of the group $ \jkf/\ohs$, which is isomorphic to the group $J_K$ of
fractional ideals.

One can define a Bost-Connes type system for
$K$~\cite{cmn1,ha-pa,LLNlat} as the corner $pAp$ of the crossed product
$$
A:=C_0(\gal(\kab / K)\times_{\ohs} \akf)\rtimes J_K,
$$
where $p$ is the characteristic function of the clopen set $\gal(\kab /
K)\times_{\ohs} \hat\OO\subset \gal(\kab / K)\times_{\ohs} \akf$. We
denote the algebra $pAp$ by $C^*_r(J_K\boxtimes(\gal(\kab /
K)\times_{\ohs} \hat\OO))$. The dynamics is defined by
$$
\sigma_t(fu_g)=N(g)^{it}fu_g\ \ \text{for}\ \ f\in C_0(\gal(\kab /
K)\times_{\ohs} \akf)\ \ \text{and}\ \ g\in J_K,
$$
where $u_g$ denotes the element of the multiplier algebra of the
crossed product corresponding to $g$, and $N\colon J_K \to (0,+\infty)$
is the absolute norm.

In~\cite{LLNlat} we showed that for every $\beta\in(0,1]$ there exists
a unique KMS$_\beta$-state $\varphi_\beta$ of the system. It is defined
by the measure $\mu_\beta$ on $\gal(\kab / K)\times_{\ohs} \akf$ which
is the push-forward of the product measure $\mu_\gal\times\prod_{v\in
V_{K,f}}\mu_{\beta,v}$ on $\gal(\kab/K)\times\akf$, where $\mu_\gal$ is
the normalized Haar measure on $\gal(\kab/K)$, and $\mu_{\beta,v}$ is
the unique measure on $K_v$ such that $\mu_{\beta,v}(\OO_v)=1$ and
$$
\mu_{\beta,v}(g
Z)=\|g\|_v^{\beta}\mu_{\beta,v}(Z)\ \ \text{for}\ \ g\in K_v^*,
$$
where $\|\cdot\|_v$ is the normalized valuation in the class $v$, so
$\|\pi\|_v=|\OO_v/\pp_v|^{-1}$ for any element~$\pi$ generating the
maximal ideal $\pp_v\subset\OO_v$.

\begin{theorem} \label{thmBC}
The KMS$_\beta$-states $\varphi_\beta$, $\beta\in(0,1]$, have type
III$_1$. In other words, for every $\beta\in(0,1]$ the action of $J_K$
on $(\gal(\kab / K)\times_{\ohs} \akf,\mu_\beta)$ is of type III$_1$.
\end{theorem}

By considering the flow of weights we can reformulate the result as
follows. Denote by $K^*_+\subset K^*$ the subgroup of totally positive
elements, that is, elements $k\in K^*$ such that $\alpha(k)>0$ for any
real embedding $\alpha\colon K\hookrightarrow \R$. Denote by
$\mu_{\beta,f}$ the measure $\prod_{v\in V_{K,f}}\mu_{\beta,v}$ on
$\akf$. Observe that for $\beta=1$ we get a Haar measure on the
additive group $\akf$.

\begin{corollary} \label{corBC}
For every $\beta\in(0,1]$, the action of $K^*_+$ on
$(\R_+\times\akf,\lambda_\infty\times\mu_{\beta,f})$ defined by
$k(t,x)=(N(k)t,kx)$, is ergodic.
\end{corollary}

\bp For any $g\in J_K$ we have $dg\mu_\beta/d\mu_\beta=N(g)^{\beta}$.
Since the action of $J_K$ on $(\gal(\kab / K)\times_{\ohs}
\akf,\mu_\beta)$ is of type III$_1$, it follows that the action of
$J_K$ on $(\R_+\times(\gal(\kab / K)\times_{\ohs}
\akf),\lambda_\infty\times\mu_\beta)$ defined by $g(t,x)=(N(g)t,gx) $
is ergodic. In other words, the action of $\akf^*$ on
$$(\R_+\times\gal(\kab /
K)\times\akf,\lambda_\infty\times\mu_\gal\times\mu_{\beta,f})$$ defined
by $g(t,x,y)=(N(g)t,s(g)^{-1}x,gy)$, is ergodic.

The Artin map $s\colon\ak^*\to\gal(\kab/K)$ is surjective with kernel
$\overline{K^*(K_\infty^*)^o}$, where $(K_\infty^*)^o$ is the connected
component of the identity in $K_\infty^*$. Since
$\ak^*=K_\infty^*\times\akf^*$, it follows that the action of
$K^*\times\akf^*$ on
$$
(\R_+\times(K^*_\infty/(K^*_\infty)^o)\times\akf^*\times\akf,
\lambda_\infty\times\nu\times\mu\times\mu_{\beta,f}),
$$
defined by $(k,g)(t,x,y,z)=(N(g)t,k^{-1}x,k^{-1}g^{-1}y,gz)$, is
ergodic, where $\nu$ and $\mu$ are Haar measures on
$K^*_\infty/(K^*_\infty)^o$ and $\akf^*$, respectively. In other words,
the action of $K^*$ on
$$
(\R_+\times(K^*_\infty/(K^*_\infty)^o)\times\akf,
\lambda_\infty\times\nu\times\mu_{\beta,f}),
$$
defined by $k(t,x,z)=(N(k)^{-1}t,k^{-1}x,k^{-1}z)$, is ergodic. Since
the group $K^*_\infty/(K^*_\infty)^o$ is finite and the homomorphism
$K^*\to K^*_\infty/(K^*_\infty)^o$ is surjective with kernel $K^*_+$,
this is equivalent to ergodicity of the action of $K^*_+$ on
$(\R_+\times\akf, \lambda_\infty\times\mu_{\beta,f})$. \ep

Observe that we can equivalently say that the action of $K^*_+$ on
$(\akf,\mu_{\beta,f})$ is ergodic of type III$_1$.

\medskip

The action of $J_K$ on $(\gal(\kab / K)\times_{\ohs} \akf,\mu_\beta)$
is ergodic by the proof of~\cite[Theorem~2.1]{LLNlat}. The computation
of the ratio set will be based on the following lemma.

\begin{lemma} \label{ldistrib}
For any $\beta\in(0,1]$, $\lambda>1$ and $\eps>0$ there exists a set
$\{\pp_n,\qq_n\}_{n\ge1}$ of prime ideals in~$\OO$ such that $$
\left|\frac{N(\qq_n)^\beta}{N(\pp_n)^\beta}-\lambda\right|<\eps \ \
\text{for all}\ \ n\ge1,\ \ \text{and}\ \
\sum^\infty_{n=1}N(\qq_n)^{-\beta}=\infty.$$
\end{lemma}

\bp The proof is the same as in~\cite{bla,boc-zah} for $K=\Q$, the only
difference is that instead of the prime number theorem one has to use
the prime ideal theorem.

It suffices to consider the case $\beta=1$. For $x>0$ denote by
$\pi(x)$ the number of prime ideals $\pp$ in~$\OO$ with $N(\pp)\le x$.
Then
$$
\pi(x)\sim\frac{x}{\log x}\ \ \text{as}\ \ x\to\infty,
$$
see e.g. \cite[Theorem 1.3]{LO}, where one can also find an estimate of
the error term.

Choose $\delta>0$ such that $1+\delta<\lambda$ and
$\lambda\delta<\eps$. Since
$$
\frac{(1+\delta)x}{\log((1+\delta)x)}-\frac{x}{\log x}\sim \frac{\delta x}{\log x},
$$
we get
\begin{equation}\label{easymp}
\pi((1+\delta)x)-\pi(x)\sim\frac{\delta x}{\log x}.
\end{equation}
In other words, if we put $B(x)=\{\pp\mid x<N(\pp)\le (1+\delta)x\}$,
then $\displaystyle |B(x)|\sim\frac{\delta x}{\log x}$. In particular,
there exists $x_0>0$ such that $|B(\lambda x)|>|B(x)|$ for all $x\ge
x_0$. Let $\pp_1,\pp_2,\dots$ be an enumeration of the set
$\cup_{m\ge0}B(\lambda^{2m}x_0)$ such that $N(\pp_1)\le
N(\pp_2)\le\dots$. For each $m\ge0$ choose a subset $C_m\subset
B(\lambda^{2m+1}x_0)$ such that $|C_m|=|B(\lambda^{2m}x_0)|$. Let
$\qq_1,\qq_2,\dots$ be an enumeration of the set $\cup_{m\ge0}C_m$ such
that $N(\qq_1)\le N(\qq_2)\le\dots$. Then, for every $n\ge1$, if
$\pp_n\in B(\lambda^{2m}x_0)$ then $\qq_n\in B(\lambda^{2m+1}x_0)$, so
that
$$
\lambda-\eps<\frac{\lambda}{1+\delta}<\frac{N(\qq_n)}{N(\pp_n)}
<\lambda(1+\delta)<\lambda+\eps.
$$
By \eqref{easymp} we also have
$$
\sum_{n=1}^\infty N(\qq_n)^{-1}\ge\sum_{m=0}^\infty
\frac{|B(\lambda^{2m}x_0)|}{(1+\delta)\lambda^{2m+1}x_0}=\infty.
$$\ep

\bp[Proof of Theorem~\ref{thmBC}] Since $\gal(\kab/K)$ is compact and
totally disconnected, it suffices to show that the action of $J_K$ on
$$
((\gal(\kab /
K)\times_{\ohs} \akf)/\gal(\kab/K),\mu_\beta)=(\akf/\ohs,\mu_{\beta,f})
$$
is of type III$_1$, see the proof of the main theorem in~\cite{nes}, as
well as~\cite[Proposition~4.6]{LLNcm} for a more general statement.

The measure space $(\hat\OO/\ohs,\mu_{\beta,f})$ can be identified with
$\prod_{v\in V_{K,f}}(\Z_+,\nu_{\beta,v})$, where $\nu_{\beta,v}$ is
the measure defined by
$$
\nu_{\beta,v}(n)=N(\pp_v)^{-n\beta}(1-N(\pp_v)^{-\beta}) \ \ \text{for}\
\ n\ge0,
$$
and the equivalence relation $\rr$ induced on $\hat\OO/\ohs$ by the
action of $J_K$ on $\akf/\ohs$ is exactly the equivalence relation
considered in our discussion of the asymptotic ratio set.

To compute $r_\infty(\rr)$, fix $\lambda>1$ and $\eps>0$. Let
$\{\pp_n,\qq_n\}_{n\ge1}$ be the set of prime ideals given by
Lemma~\ref{ldistrib}. Let $v_n$ and $w_n$ be the places corresponding
to $\pp_n$ and $\qq_n$, respectively. Then we define the sets
$I_n\subset V_{K,f}$ and $K_n,L_n\subset\prod_{v\in I_n}\Z_+$ required
by the definition of the asymptotic ratio set by
$$
I_n=\{v_n,w_n\},\ \ K_n=\{(0,1)\}, \ \ L_n=\{(1,0)\},
$$
and denote by $\varphi_n\colon K_n\to L_n$ the unique bijection. For
$a=(0,1)\in K_n$ and $b=(1,0)\in L_n$ we have
$$
\frac{\mu_{\beta,f}(Z(b))}{\mu_{\beta,f}(Z(a))}
=\frac{\nu_{\beta,v_n}(1)\nu_{\beta,w_n}(0)}{\nu_{\beta,v_n}(0)\nu_{\beta,w_n}(1)}
=\frac{N(\pp_n)^{-\beta}(1-N(\pp_n)^{-\beta})(1-N(\qq_n)^{-\beta})}
{N(\qq_n)^{-\beta}(1-N(\pp_n)^{-\beta})(1-N(\qq_n)^{-\beta})},
$$
which for large $n$ is arbitrarily close to
$N(\qq_n)^\beta/N(\pp_n)^\beta$, which in turn is close to $\lambda$ up
to $\eps$. We also have
$$
\sum_{n=1}^\infty\sum_{a\in K_n}\mu_{\beta,f}(Z(a))
=\sum^\infty_{n=1}N(\qq_n)^{-\beta}(1-N(\pp_n)^{-\beta})(1-N(\qq_n)^{-\beta})
=\infty.
$$
It follows that $\lambda\in r_\infty(\rr)$. Since this is true for all
$\lambda>1$, we conclude that the action is of type~III$_1$. \ep

\begin{remark}
In view of Corollary~\ref{corBC} it is natural to ask whether the
action of $K^*$ on
$(K_\infty\times\akf,\bar\lambda_\infty\times\mu_{\beta,f})$ is
ergodic, where $\bar\lambda_\infty$ is a Haar measure on $K_\infty\cong
\R^{[K\colon\Q]}$; see also Remark~\ref{rgln}(ii) below for a more
general question. Assume for simplicity that $K$ is an imaginary
quadratic field of class number one. Then one can try to prove that the
action of $K^*$ on $\T\times\akf$ given by $k(z,x)=(k|k|^{-1}z,kx)$, is
ergodic, e.g. by adapting the strategy in \cite{nes,LLNlat}, and then
compute the ratio set of this action. However, for the latter one would
need an information about the distribution of the angles of prime
ideals, that is, of the values of the homomorphism $J_K\to \T/\OO^*$,
$(k)\to k|k|^{-1}$. We are not aware of any result of this sort.
\end{remark}

\bigskip
\section{The Connes-Marcolli $\gl$-system}

Let $G$ be a discrete group and $\Gamma$ be a subgroup of $G$. Recall
that $(G,\Gamma)$ is called a Hecke pair if every double coset of
$\Gamma$ contains finitely many right cosets of $\Gamma$, so that
$$
R_\Gamma(g):=|\Gamma\backslash\Gamma g\Gamma|<\infty\ \ \text{for
any}\ \ g\in G.
$$
Then the space $\hecke{G}{\Gamma}$ of finitely supported functions on
$\Gamma\backslash G/\Gamma$ is a $*$-algebra with  product
$$
(f_1*f_2)(g)=\sum_{h\in\Gamma\backslash G}f_1(gh^{-1})f_2(h)
$$
and involution $f^*(g)=\overline{f(g^{-1})}$. Denote by $[g]\in
\hecke{G}{\Gamma}$ the characteristic function of the double coset
$\Gamma g\Gamma$.

If $G$ acts on a space $X$ then every element $g\in G$ defines a Hecke
operator $T_g$ acting on functions on $\Gamma\backslash X$, which we
also consider as $\Gamma$-invariant functions on $X$:
$$
(T_gf)(x)=\frac{1}{R_\Gamma(g)}\sum_{h\in\Gamma\backslash\Gamma
g\Gamma}f(hx).
$$
Then $[g^{-1}]\mapsto R_\Gamma(g)T_g$ is a representation of the Hecke
algebra $\hecke G \Gamma$ on the space of $\Gamma$-invariant functions.

If $X$ is locally compact and the action of $\Gamma$ on $X$ is proper,
one can define a C$^*$-algebra $C^*_r(\fib{G}{\Gamma}{X})$ which can be
thought of as a crossed product of $C_0(\Gamma\backslash X)$ by $\hecke
G \Gamma$, see~\cite{con-mar,LLNcm}. It is a completion of the algebra
$C_c(\fib{G}{\Gamma}{X})$ of continuous compactly supported functions
on $\fib{G}{\Gamma}{X}$ with convolution product
\begin{equation*} \label{econv}
(f_1*f_2)(g,x)=\sum_{h\in\Gamma\backslash G}f_1(gh^{-1},hx)
f_2(h,x)
\end{equation*}
and involution $f^*(g,x)=\overline{f(g^{-1},gx)}$. If the action of
$\Gamma$ is free then $\fib{G}{\Gamma}{X}$ is a groupoid and
$C^*_r(\fib{G}{\Gamma}{X})$ is the usual groupoid C$^*$-algebra.

\medskip

Consider now the Hecke pair $(\glq,\slz)$, where $\glq$ is the group of
rational matrices with positive determinant. The group $\glq$ acts by
multiplication on $\mtwo(\Q_p)$ for every prime~$p$. It also acts by
M\"obius transformations on the upper half-plane $\HH$. The
$\gl$-system of Connes and Marcolli~\cite{con-mar} is the corner of the
C$^*$-algebra
$$C^*_r(\fib{\glq}{\slz}{(\HH\times\mtwo(\af))})$$ defined by the
projection corresponding to the subspace
$\HH\times\mtwo(\hat\Z)\subset\HH\times\ma$, where $\af=\A_{\Q,f}$. We denote this algebra by
$\red{\fibb{\glq}{\slz}{(\HH\times\mtwo(\hat\Z))}}$. The dynamics on it
is defined by
$$
\sigma_t(f)(g,x)=\det(g)^{it}f(g,x).
$$

In~\cite{LLNcm} we showed that for every $\beta\in(1,2]$ there exists a
unique KMS$_\beta$-state $\varphi_\beta$ on the $\gl$-system. It is
defined by the product-measure $\mu_\HH\times\prod_p\mu_{\beta,p}$ on
$\HH\times\ma$, where $\mu_\HH$ is the unique $\glp(\R)$-invariant
measure on $\HH$ such that $\mu_\HH(\slz\backslash\HH)=2$, and
$\mu_{\beta,p}$ is the unique measure on $\mtwo(\Q_p)$ such that
$\mu_{\beta,p}(\mtwo(\Z_p))=1$ and
$$
\mu_{\beta,p}(gZ)=|\det(g)|_p^{\beta}\mu_{\beta,p}(Z)\ \ \text{for}\ \
g\in\gl(\Q_p).
$$
Denote the measure $\prod_p\mu_{\beta,p}$ by $\mu_{\beta,f}$. Observe
that $\mu_{2,f}$ is the Haar measure of the additive group $\ma$
normalized so that $\mu_{2,f}(\mtwo(\hat\Z))=1$.

The definition of the $\gl$-system required a new type of crossed
product construction because of non-freeness of the action of $\slz$ on
$\HH\times\ma$. However, the set of points with nontrivial stabilizers
is $\HH\times\{0\}$, which has measure zero with respect to
$\mu_\HH\times\mu_{\beta,f}$. As a result the von Neumann algebra
generated by the $\gl$-system in the GNS-representation
of~$\varphi_\beta$ is much easier to describe. It is the reduction of
the von Neumann algebra crossed product
$L^\infty(\HH\times\ma,\mu_\HH\times\mu_{\beta,f})\rtimes\glq$ by the
projection defined by a fundamental domain of the action of $\slz$ on
$\HH\times(\mtwo(\hat\Z)\setminus\{0\})$. Therefore to compute the type
of the algebra we have to compute the type of the action of $\glq$ on
$\HH\times\ma$.

It is natural to consider a slightly more general problem. Namely,
replace $\HH=\pgl/\pso$ by $\pgll$ and $\glq$ by $\gl(\Q)$. Denote by
$\mu_\infty$ the Haar measure of $\pgll$ normalized so that
$\mu_\infty(\gl(\Z)\backslash\pgll)=2$. Put
$\mu_\beta=\mu_\infty\times\mu_{\beta,f}$. The action of $\gl(\Q)$ on
$(\pgll\times\ma,\mu_\beta)$ is ergodic by~\cite[Corollary~4.7]{LLNcm}.

\begin{theorem} \label{thmgl2}
For every $\beta\in(1,2]$, the action of $\gl(\Q)$ on
$(\pgll\times\ma,\mu_\beta)$ has type~III$_1$. In particular, the
KMS$_\beta$-states $\varphi_\beta$, $\beta\in(1,2]$, of the
Connes-Marcolli $\gl$-system have type~III$_1$.
\end{theorem}

As we already remarked in~\cite{LLNcm}, the flows of weights of the
above actions are easy to describe, and then the result takes the
following essentially equivalent form. Denote by $\lambda_\infty$ the
usual Lebesgue measure on $\mtwo(\R)\cong\R^4$, and put
$\lambda_\beta=\lambda_\infty\times\mu_{\beta,f}$. For $\beta=2$ we get
a Haar measure on the additive group $\mtwo(\A)=\mtwo(\R)\times\ma$,
where $\A=\A_\Q$.

\begin{corollary}
For every $\beta\in(1,2]$, the action of $\gl(\Q)$ on
$(\mtwo(\A),\lambda_\beta)$ is ergodic.
\end{corollary}

\bp If the action of $\gl(\Q)$ on $(\pgll\times\ma,\mu_\beta)$ is of
type III$_1$ then clearly also the action of $\glq$ on
$(\pgl\times\ma,\mu_\beta)$ is of type III$_1$. As we discussed
in~\cite[Remark~4.9]{LLNcm}, using the isomorphism
$\glp(\R)/\{\pm1\}\cong\R^*_+\times\pgl$ we can identify the underlying
space of the flow of weights of this action with the quotient of
$((\glp(\R)/\{\pm1\})\times\ma,\lambda_\beta)$ by the diagonal action
of $\glq$. Therefore this diagonal action is ergodic. Since $\glp(\R)$
is connected and $\{\pm1\}$ is finite, by~\cite[Proposition~4.6]{LLNcm}
we conclude that the action of $\glq$ on
$(\glp(\R)\times\ma,\lambda_\beta)$ is ergodic. But then the action of
$\gl(\Q)$ on $(\gl(\R)\times\ma,\lambda_\beta)$ is also ergodic. \ep

To simplify notation from now on we write $G$ for $\gl(\Q)$, $\Gamma$
for $\gl(\Z)$ and $X$ for $\pgll\times\ma$.

Recall, see e.g.~\cite{kri}, that the group $G$ is generated by
$\Gamma$ and the matrices $\diag 1 p$, $p\in\primes$, where $\primes$
is the set of prime numbers. We have
\begin{equation} \label{egenp}
\{m\in\mz: |\det(m)|=p\}=\Gamma \diag{1}{p}\Gamma\ \ \text{and}\ \
R_\Gamma\left(\diag{1}{p}\right)=p+1.
\end{equation}

Recall also, see~\cite[Section 3]{LLNcm}, that
\begin{equation} \label{enormalp}
\mu_{\beta,p}(\gl(\Z_p))=(1-p^{-\beta})(1-p^{-\beta+1}),
\end{equation}
Using the scaling property of $\mu_{\beta,p}$ and \eqref{egenp} we then
conclude that
\begin{equation}\label{enormalp2}
\mu_{\beta,p}(\{m\in\mtwo(\Z_p):|\det(m)|_p=p^{-1}\})=p^{-\beta}(p+1)
(1-p^{-\beta})(1-p^{-\beta+1}).
\end{equation}

\medskip

One of the key ingredients of the proof of the theorem will be the
following version of equidistribution of Hecke points. I am indebted to
Hee Oh for explaining me how to put it in the general setup of
equidistibution of Hecke points for reductive groups.

Denote by $G_p$ the group generated by $\Gamma$ and $\diag 1 p$.
Equivalently, $G_p$ is the group $\gl(\Z[p^{-1}])$. For a nonempty
subset $F$ of primes denote by $G_F$ the group generated by $G_p$ for
all $p\in F$.
For finite~$F$ denote by~$X_F$ the space $\pgll\times\mtwo(\Q_F)$,
where $\Q_F=\prod_{p\in F}\Q_p$, and by $\bar\mu_{\beta,F}$ the measure
$\mu_\infty\times\prod_{p\in F}\mu_{\beta,p}$. Consider also the
measure~$\bar\nu_{\beta,F}$ on $\Gamma\backslash X_F$ defined by
$\bar\mu_{\beta,F}$. We shall write~$\Z_F$ for~$\prod_{p\in F}\Z_p$.

\begin{lemma} \label{lequi}
Let $F$ be a finite set of primes, $r\in\gl(\Q_F)$. Assume $f$ is a
compactly supported continuous right $\gl(\Z_F)$-invariant function on
$Z=\Gamma\backslash(\pgll\times\gl(\Z_F)r\gl(\Z_F))\subset\Gamma\backslash
X_F$. Then for any $\eps>0$ and any compact subset $C$ of $Z$ there
exists $M>0$ such that if $g\in G_{F^c}$ ($F^c=\primes\setminus F$) and
$R_\Gamma(g)\ge M$ then
$$
\left|T_gf(x)-\bar\nu_{\beta,F}(Z)^{-1}\int_Zf\,d\bar\nu_{\beta,F}\right|<\eps \ \
\text{for all}\ \ x\in C.
$$
\end{lemma}

\bp The $\gl(\Z_F)$-space $\gl(\Z_F)r\gl(\Z_F)/\gl(\Z_F)$ can be
identified with $\gl(\Z_F)/H$, where $H=r\gl(\Z_F)r^{-1}\cap\gl(\Z_F)$.
Therefore $f$ can be considered as a function on
$\Gamma\backslash(\pgll\times\gl(\Z_F))/H$. Consider the compact open
subgroup $U=H\times\prod_{p\in F^c}\gl(\Z_p)$ of $\gl(\hat\Z)$. By
considering $\gl(\Z_F)$ as the subgroup of $\gl(\hat\Z)$ consisting of
elements with coordinates $1$ for $p\in F^c$, we get a homeomorphism
$$
\Gamma\backslash(\pgll\times\gl(\Z_F))/H\to
G\backslash(\pgll\times\gl(\af))/U,
$$
since $G\cap\gl(\hat \Z)=\Gamma$ and $\gl(\af)=G\gl(\hat\Z)$.
Furthermore, we have $\gl(\af)=GU$. To see this recall that
$\sltwo(\Q)$ is dense in $\sltwo(\af)$ by the strong approximation
theorem. It follows that $GU$ contains $\sltwo(\af)$, and in particular
$\sltwo(\hat\Z)$. Since $\gl(\af)=G\gl(\hat\Z)$, it is therefore enough
to check that $\sltwo(\hat\Z)U=\gl(\hat\Z)$, that is,
$\sltwo(\Z_F)H=\gl(\Z_F)$. For this we have to show that the
determinant map $\det\colon H\to\Z_F^*$ is surjective. Since every
double coset of $\gl(\Z_F)$ contains a diagonal matrix, without loss of
generality we may assume that $r$ is diagonal. But then $H$ contains
all the diagonal matrices of $\gl(\Z_F)$, and surjectivity of the
determinant is immediate.

We can then proceed as in \cite{cou}, see Remark (3) following
\cite[Theorem 1.7]{cou}, as well as Sections 2 and 3 in \cite{go}. \ep

For a $\Gamma$-invariant measurable subset $A$ of $\pgll\times\ma$,
denote by $m(A)$ the operator of multiplication by the characteristic
function of $A$ on $L^2(\Gamma\backslash X,d\nu_\beta)$, where
$\nu_\beta$ is the measure on $\Gamma\backslash X$ defined by
$\mu_\beta$.

\begin{lemma} \label{lbasicsets}
For a prime $p$, consider the sets
\begin{align*}
A_p&=\{x\in \pgll\times\mtwo(\hat\Z): x_p\in\gl(\Z_p)\},\\
B_p&=\{x\in
\pgll\times\mtwo(\hat\Z): |\det(x_p)|_p=p^{-1}\},
\end{align*}
and the element $g=\diag{1}{p}$. Then for the
operator $m(A_p)T_gm(B_p)=m(A_p)T_g$ on $L^2(\Gamma\backslash X,d\nu_\beta)$ we have
$$
\|m(A_p)T_gm(B_p)\|=p^{\beta/2}(p+1)^{-1/2}=
\nu_\beta(\Gamma\backslash A_p)^{1/2}\nu_\beta(\Gamma\backslash B_p)^{-1/2}.
$$
\end{lemma}

\bp Since $B_p=\Gamma gA_p$, $R_\Gamma(g)=p+1$ and $|T_g(f)|^2\le
T_g(|f|^2)$ pointwise by convexity, this follows from
\cite[Lemma 2.7]{LLNcm}, but we will sketch a proof for the reader's
convenience.

Fix representatives $h_1,\dots,h_{p+1}$ of $\Gamma\backslash\Gamma
g\Gamma$. Choose a fundamental domain $C$ for the action of $\Gamma$ on
$A_p$. Using that the action of $\Gamma$ on $A_p$ is free and that
$G_p\cap\gl(\Z_p)=\Gamma$ one can easily check that the sets $\Gamma
h_iC$ are mutually disjoint and the factor-map $p\colon X\to
\Gamma\backslash X$ is injective on the sets~$h_iC$. Consider the
operators $S_i$ defined by
$$
(S_if)(p(x))=\begin{cases}f(p(h_ix)), &\text{if }x\in C,\\
0, &\text{if }x\notin A_p.\end{cases}
$$
Then $p^{-\beta/2}S_i$ is a partial isometry with initial space
$L^2(p(h_iC),d\nu_\beta)$ and range  $L^2(p(C),d\nu_\beta)$, and
$m(A_p)T_gm(B_p)=(S_1+\dots+S_{p+1})/(p+1)$. Since the spaces
$L^2(p(h_iC),d\nu_\beta)$ are mutually orthogonal, we have
$$
\|p^{-\beta/2}S_1+\dots+p^{-\beta/2}S_{p+1}\|=(p+1)^{1/2},
$$
which gives the first equality in the statement. The second equality
follows from \eqref{enormalp} and \eqref{enormalp2}.\ep

For $x\in X$ denote by $\bar x_F$ the image of $x$ under the factor-map
$X\to X_F$. For a function $f$ on $X_F$ consider the function $f_F$ on
$X$ defined by
$$
f_F(x)=\begin{cases}f(\bar x_F), &\text{if }x_p\in\mtwo(\Z_p)
\text{ for all }p\notin F,\\ 0,&\text{otherwise}.\end{cases}
$$

\begin{lemma} \label{lmaingl2}
For any $\beta\in(1,2]$ and $\lambda>1$ there exists $c>0$ such that
for any $\eps>0$, any finite set~$F$ of primes and any positive
compactly supported continuous right $\gl(\Z_F)$-invariant function~$f$
on $\Gamma\backslash(\pgll\times\mtwo(\Z_F))\subset\Gamma\backslash
X_F$ with $\int_{\Gamma\backslash X_F}fd\bar\nu_{\beta,F}=1$, there
exist a subset $\{p_n,q_n\}_{n\ge1}$ of $F^c$ and $\Gamma$-invariant
measurable subsets $X_{1n}$, $X_{2n}$, $Y_{1n}$ and $Y_{2n}$, $n\ge1$,
of $X$ such that \enu{i}
$\displaystyle\left|\frac{q_n^\beta}{p_n^\beta}-\lambda\right|<\eps$
for all $n\ge1$; \enu{ii} the sets $Y_{1n}$, $n\ge1$, as well as the
sets $Y_{2n}$, $n\ge1$, are mutually disjoint; \enu{iii}\ \
$\displaystyle\sum^\infty_{n=1}
\left(\frac{m(X_{1n})T_{g_n}m(Y_{1n})}{\|m(X_{1n})T_{g_n}m(Y_{1n})\|}
f_F,\frac{m(X_{2n})T_{h_n}m(Y_{2n})}{\|m(X_{2n})T_{h_n}m(Y_{2n})\|}
f_F\right)_{L^2(\Gamma\backslash X,d\nu_\beta)}>c$, where
$g_n=\diag{1}{p_n}$ and $h_n=\diag{1}{q_n}$.
\end{lemma}

\bp Fix $\delta\in(0,1)$. Choose representatives $r_k$, $k\ge 1$, of
the double cosets
$$
\gl(\Z_F)\backslash(\gl(\Q_F)\cap\mtwo(\Z_F))/\gl(\Z_F),
$$
and put $Z_k=\Gamma\backslash(\pgll\times \gl(\Z_F)r_k\gl(\Z_F))$. Let
$N$ be such that
$$
\sum^N_{k=1}\int_{Z_k}fd\bar\nu_{\beta,F}
>\int_{\Gamma\backslash X_F}fd\bar\nu_{\beta,F}-\delta=1-\delta.
$$
Next choose compact subsets $C_k$ of $Z_k$ such that
$$
\bar\nu_{\beta,F}(C_k)>(1-\delta)\bar\nu_{\beta,F}(Z_k)
\ \ \text{for}\ \ k=1,\dots,N.
$$
By Lemma~\ref{lequi} there exists $M$ such that for any element $g\in
G_{F^c}$ with $R_\Gamma(g)\ge M$ we have
$$
T_gf(x)\ge\frac{1-\delta}{\bar\nu_{\beta,F}(Z_k)}
\int_{Z_k}fd\bar\nu_{\beta,F}\ \
\text{for}\ \ x\in C_k,\ \ k=1,\dots,N.
$$
It follows that if we take two elements $g,h\in G_{F^c}$ with
$R_\Gamma(g),R_\Gamma(h)\ge M$, then
\begin{align*}
\int_{\Gamma\backslash X_F}T_gf\,T_hfd\bar\nu_{\beta,F}&\ge
\sum^N_{k=1}\int_{C_k}T_gf\,T_hfd\bar\nu_{\beta,F}\\
&\ge \sum^N_{k=1}\left(\frac{1-\delta}{\bar\nu_{\beta,F}(Z_k)}
\int_{Z_k}fd\bar\nu_{\beta,F}\right)^2\bar\nu_{\beta,F}(C_k)\\
&\ge (1-\delta)^3\sum^N_{k=1}\left(\frac{1}{\bar\nu_{\beta,F}(Z_k)}
\int_{Z_k}fd\bar\nu_{\beta,F}\right)^2
\bar\nu_{\beta,F}(Z_k)\\
&\ge (1-\delta)^3\left(
\sum^N_{k=1}\int_{Z_k}fd\bar\nu_{\beta,F}\right)^2,
\end{align*}
since the function $t\mapsto t^2$ is convex and
$\sum^N_{k=1}\bar\nu_{\beta,F}(Z_k)\le
\bar\nu_{\beta,F}(\Gamma\backslash (\pgll\times\mtwo(\Z_F)))=1$.
Therefore
\begin{equation} \label{eheckemix}
\int_{\Gamma\backslash X_F}T_gf\,T_hfd\bar\nu_{\beta,F}
\ge(1-\delta)^5.
\end{equation}

Using Lemma~\ref{ldistrib} choose a subset $\{p_n,q_n\}_{n\ge1}$ of
$F^c$ such that $q_n>p_n\ge M$ and $|q_n^\beta/
p_n^\beta-\lambda|<\eps$ for all~$n$, and
\begin{equation} \label{einfprimes}
\sum_n p_n^{1-\beta}=\infty.
\end{equation}
Consider the sets $A_p$ and $B_p$ from Lemma~\ref{lbasicsets} and put
$$
X_{1n}=A_{p_n}\backslash(B_{p_1}\cup\dots\cup B_{p_{n-1}}),\ \
Y_{1n}=B_{p_n}\backslash(B_{p_1}\cup\dots\cup B_{p_{n-1}}), 
$$
$$
X_{2n}=A_{q_n}\backslash(B_{q_1}\cup\dots\cup B_{q_{n-1}}),\ \
Y_{2n}=B_{q_n}\backslash(B_{q_1}\cup\dots\cup B_{q_{n-1}}). 
$$
Let $g_n=\diag{1}{p_n}$ and $h_n=\diag{1}{q_n}$.

We claim that if $g\in\Gamma g_n\Gamma$ and $x\in X_{1n}$ then $gx\in
Y_{1n}$. Indeed, we clearly have $g A_{p_n}\subset B_{p_n}$, so that
$gx\in B_{p_n}$. Furthermore, if $gx\in B_{p_k}$ for some $k<n$ then
$p_k^{-1}=|\det(gx_{p_k})|_{p_k}=|\det(x_{p_k})|_{p_k}$, since
$g_n\in\gl(\Z_{p_k})$. Therefore $x\in B_{p_k}$, which contradicts the
assumption that $x\in X_{1n}$. Hence $gx\in Y_{1n}$.

It follows that
$m(X_{1n})T_{g_n}m(Y_{1n})f_F=(T_{g_n}f)_F\31_{\Gamma\backslash
X_{1n}}$. For similar reasons,
$m(X_{2n})T_{h_n}m(Y_{2n})f_F=(T_{h_n}f)_F\31_{\Gamma\backslash
X_{2n}}$. Therefore
\begin{multline*}
(m(X_{1n})T_{g_n}m(Y_{1n})f_F,m(X_{2n})
T_{h_n}m(Y_{2n})f_F)_{L^2(\Gamma\backslash X,d\nu_\beta)}\\
=(T_{g_n}f,T_{h_n}f)_{L^2(\Gamma\backslash X_F,d\bar\nu_{\beta,F})}
\nu_\beta(\Gamma\backslash(X_{1n}\cap X_{2n}))
\ge(1-\delta)^5\nu_\beta(\Gamma\backslash(X_{1n}\cap X_{2n}))
\end{multline*}
by \eqref{eheckemix}.

By Lemma~\ref{lbasicsets} we have
$$
\|m(X_{1n})T_{g_n}m(Y_{1n})\|\le
\nu_\beta(\Gamma\backslash B_{p_n})^{-1/2}\ \ \text{and}\ \
\|m(X_{2n})T_{h_n}m(Y_{2n})\|\le \nu_\beta(\Gamma\backslash B_{q_n})^{-1/2}.
$$
If follows that
\begin{multline}\label{eright}
\sum^\infty_{n=1}
\left(\frac{m(X_{1n})T_{g_n}m(Y_{1n})}{\|m(X_{1n})T_{g_n}m(Y_{1n})\|}f_F,
\frac{m(X_{2n})T_{h_n}m(Y_{2n})}{\|m(X_{2n})T_{h_n}m(Y_{2n})\|}
f_F\right)_{L^2(\Gamma\backslash
X,d\nu_\beta)}\\
\ge(1-\delta)^5\sum^\infty_{n=1}(\nu_\beta(\Gamma\backslash B_{p_n})
\nu_\beta(\Gamma\backslash B_{q_n}))^{1/2}\nu_\beta(\Gamma\backslash(X_{1n}\cap
X_{2n})).
\end{multline}
We have
\begin{equation} \label{emeasure1}
\nu_\beta(\Gamma\backslash(X_{1n}\cap
X_{2n}))=\nu_\beta(\Gamma\backslash A_{p_n})\nu_\beta(\Gamma
\backslash A_{q_n})\prod^{n-1}_{k=1}
(1-\nu_\beta(\Gamma\backslash( B_{p_k}\cup B_{q_k}))).
\end{equation}
We may assume that $M$ is so large that
\begin{equation} \label{emeasure2}
\nu_\beta(\Gamma\backslash A_{p_n})\nu_\beta(\Gamma
\backslash A_{q_n})
=(1-p_n^{-\beta})(1-p_n^{-\beta+1})(1-q_n^{-\beta})(1-q_n^{-\beta+1})
>1-\delta,
\end{equation}
see~\eqref{enormalp}. Since
$$
\nu_\beta(B_p)=p^{-\beta}(p+1)(1-p^{-\beta})(1-p^{-\beta+1})\sim p^{1-\beta}
$$
by~\eqref{enormalp2}, we may also assume that $M$ is so large and the
ratios $q_n^\beta/p_n^\beta$ are so close to $\lambda$ that
\begin{align}
(\nu_\beta(\Gamma\backslash B_{p_n})\nu_\beta(\Gamma\backslash B_{q_n}))^{1/2}
&>(1-\delta)c_0\big(\nu_\beta(\Gamma\backslash
B_{p_n})+ \nu_\beta(\Gamma\backslash
B_{q_n})-\nu_\beta(\Gamma\backslash
 (B_{p_n})\nu_\beta(\Gamma\backslash B_{q_n})\big)\notag\\
 &=(1-\delta)c_0\nu_\beta(\Gamma\backslash (B_{p_n}\cup B_{q_n})),\label{emeasure3}
\end{align}
where
$$
c_0=\frac{\lambda^{(1-\beta)/2\beta}}{1+\lambda^{(1-\beta)/\beta}}.
$$

Combining \eqref{emeasure1}-\eqref{emeasure3} we conclude that the
right hand side of \eqref{eright} is not smaller than
$$
(1-\delta)^7c_0\sum^\infty_{n=1}
\nu_\beta(\Gamma\backslash
 (B_{p_n}\cup B_{q_n}))
\prod^{n-1}_{k=1}
(1-\nu_\beta(\Gamma\backslash( B_{p_k}\cup B_{q_k})))
=(1-\delta)^7c_0,
$$
because
$$
\sum^\infty_{n=1}\nu_\beta(\Gamma\backslash (B_{p_n}\cup B_{q_n}))
\ge \sum^\infty_{n=1}\nu_\beta(\Gamma\backslash B_{p_n})
=\infty
$$
by~\eqref{einfprimes}. Since $\delta$ can be arbitrarily small, we see
that we can take any $c<c_0$.
 \ep

\bp[Proof of Theorem~\ref{thmgl2}] Similarly to the proof of
Theorem~\ref{thmBC}, since $\gl(\hat\Z)$ is compact and totally
disconnected, it suffices to show that the action of $\gl(\Q)$ on
$(\pgll\times\mtwo(\af)/\gl(\hat\Z),\mu_\beta)$ is of type III$_1$. In
other words, in computing the ratio set it suffices to consider right
$\gl(\hat\Z)$-invariant sets.

Let $\lambda>1$, $\eps>0$ and $Y$ be a measurable right
$\gl(\hat\Z)$-invariant subset of $X$ such that $\mu_\beta(Y)>0$. Let
$c>0$ be given by Lemma~\ref{lmaingl2}. There exists $g_0\in\gl(\Q)$
such that the intersection $Y_0:=g_0Y\cap(\pgll\times\mtwo(\hat \Z))$
has positive measure. Let $\varphi$ be the function
$\nu_\beta(\Gamma\backslash \Gamma Y_0)^{-1}\31_{\Gamma\backslash
\Gamma Y_0}$ on~$\Gamma\backslash X$. We can find a finite set $F$ of
primes and a positive compactly supported continuous right
$\gl(\Z_F)$-invariant function $f$ on
$\Gamma\backslash(\pgll\times\mtwo(\Z_F))$ such that
$$
\int_{\Gamma\backslash X_F}f\,d\bar\nu_{\beta,F}=1\ \ \text{and}\ \
\|\varphi-f_F\|_2(\|\varphi\|_2+\|f_F\|_2)<c.
$$
Let $p_n,q_n\in F^c$, $X_{1n},X_{2n},Y_{1n},Y_{2n}\subset X$ and
$g_n,h_n\in G$ be given by Lemma~\ref{lmaingl2}.

\smallskip

Denote by $T'_n$ and $T''_n$ the contractions $\displaystyle
\frac{m(X_{1n})T_{g_n}m(Y_{1n})}{\|m(X_{1n})T_{g_n}m(Y_{1n})\|}$ and
$\displaystyle\frac{m(X_{2n})T_{h_n}m(Y_{2n})}{\|m(X_{2n})T_{h_n}m(Y_{2n})\|}$,
respectively, and by $e'_n$ and $e''_n$ the projections $m(Y_{1n})$ and
$m(Y_{2n})$. We have
\begin{align*}
(T'_n\varphi,T''_n\varphi)
&\ge (T'_nf_F,T''_nf_F)-\|T'_n\varphi-T'_nf_F\|_2\|T''_n\varphi\|_2
-\|T''_n\varphi-T''_nf_F\|_2\|T'_nf_F\|_2\\
&\ge (T'_nf_F,T''_nf_F)-\|e'_n(\varphi-f_F)\|_2\|e''_n\varphi\|_2
-\|e''_n(\varphi-f_F)\|_2\|e'_nf_F\|_2.
\end{align*}
Since the projections $e'_n$, as well as the projections $e''_n$, are
mutually orthogonal, we have
$$
\sum_n\|e'_n(\varphi-f_F)\|_2\|e''_n\varphi\|_2\le
\left(\sum_n\|e'_n(\varphi-f_F)\|_2^2\right)^{1/2}
\left(\sum_n\|e''_n\varphi\|_2^2\right)^{1/2}\le
\|\varphi-f_F\|_2\|\varphi\|_2,
$$
and similarly
$$
\sum_n\|e''_n(\varphi-f_F)\|_2\|e'_nf_F\|_2\le \|\varphi-f_F\|_2\|f_F\|_2.
$$
It follows that
$$
\sum_n(T'_n\varphi,T''_n\varphi)
\ge \sum_n(T'_nf_F,T''_nf_F)-\|\varphi-f_F\|_2(\|\varphi\|_2+\|f_F\|_2)>c-c=0.
$$
Hence there exists $n$ such that $(T'_n\varphi,T''_n\varphi)>0$. Then
$$
\int_{\Gamma\backslash X}T_{g_n}\varphi\,T_{h_n}\varphi\, d\nu_\beta>0,
$$
which means that the set $\Gamma g_n^{-1}\Gamma Y_0\cap \Gamma
h_n^{-1}\Gamma Y_0$ has positive measure. It follows that there exist
$g\in\Gamma g_n\Gamma$ and $h\in\Gamma h_n\Gamma$ such that
$g^{-1}Y_0\cap h^{-1}Y_0$ has positive measure, and hence
$g_0^{-1}hg^{-1}g_0Y\cap Y$ has positive measure. Since
$$
\frac{d(g_0^{-1}hg^{-1}g_0\mu_\beta)}{d\mu_\beta}
=|\det(g_0^{-1}hg^{-1}g_0)|^\beta
=\det(h_ng^{-1}_n)^\beta
=\frac{q_n^\beta}{p_n^\beta},
$$
and $|q_n^\beta/ p_n^\beta-\lambda|<\eps$, we conclude that $\lambda$
belongs to the ratio set of the action of $\gl(\Q)$ on
$(\pgll\times\mtwo(\af)/\gl(\hat\Z),\mu_\beta)$. Since this is true for
all $\lambda>1$, the action is of type III$_1$. \ep

We finish our discussion of the $\gl$-system with the following simple
observation.

\begin{proposition}
For every $\beta\in(1,2]$, the action of $\gl(\Q)$ on $(\pgll\times\ma,
\mu_\beta)$ is amenable. Therefore  the von Neumann algebra generated
by the $\gl$-system in the GNS-represen\-tation of~$\varphi_\beta$,
$\beta\in(1,2]$, is the injective factor of type III$_1$.
\end{proposition}

\bp If $F$ is a finite set of primes then $G_F$ is a discrete subgroup
of $\pgll\times\gl(\Q_F)$, since $G_F\cap\gl(\Z_F)=\Gamma$ and the
homomorphism $\Gamma\to\pgll$ has discrete image and finite kernel.
Since $\pgll\times\gl(\Q_F)$ is a subset of $X_F$ of full measure, it
follows that $L^\infty(X_F,\bar\mu_{\beta,F})\rtimes G_F$ is a type I
von Neumann algebra. Since the algebra $L^\infty(X,\mu_\beta)\rtimes G$
is the closure of an increasing union of algebras of the form
$L^\infty(X_F,\bar\mu_{\beta,F})\rtimes G_F$, it is injective, and
therefore the action of $G$ on $(X,\mu_\beta)$ is
amenable~\cite{zim_hyp}. \ep

\begin{remark} \label{rgln}
\mbox{\ } \enu{i} Since the C$^*$-algebra $A$ of the $\gl$-system is a
subalgebra of the crossed product of
$C_0(\slz\backslash(\HH\times\ma))$ by the Hecke algebra
$\hecke{\glq}{\slz}$, which is abelian, one might expect that not only
the von Neumman algebras~$\pi_{\varphi_\beta}(A)''$ are injective, but
that $A$ is nuclear. To see that this is not the case, consider the
state $\varphi_0$ on $A$ defined by the measure $\frac{1}{2}\mu_\HH$,
considered as a measure supported on
$\HH\times\{0\}\subset\HH\times\ma$. The algebra $\pi_{\varphi_0}(A)''$
is a reduction of the crossed product
$L^\infty(\HH,\mu_\HH)\rtimes(\glq/\{\pm1\})$. Therefore if $A$ were
nuclear, the action of $\glq/\{\pm1\}$ on $(\HH,\mu_\HH)$ would be
amenable, which would contradict~\cite[Corollary~1.2]{zim_dense}.
\enu{ii} It is apparently straightforward to extend the above results
to $\gln$ (with the interval~$(1,2]$ replaced by $(n-1,n]$). One can
however formulate a more general problem. Let $K$ be a number field,
$M$ a finite dimensional central simple $K$-algebra, $G$ the group of
invertible elements in~$M$.
Is the action of $G(K)$ on $M(\ak)$ with its Haar measure ergodic?
\end{remark}

\end{document}